\documentclass{article}
\usepackage{amsmath}
\usepackage{amssymb}
\usepackage{graphics}
\usepackage{color}

\def\inv{\preceq}
\def\Sub{\mbox{\rm Sub}}
\def\Nat{\mathbb{N}}
\def\A{\mathcal{A}}
\def\M{\mathcal{M}}

\newtheorem{theorem}{Theorem}[section]
\newtheorem{corollary}[theorem]{Corollary}
\newtheorem{lemma}[theorem]{Lemma}
\newtheorem{proposition}[theorem]{Proposition}

\newtheorem{ddefinition}[theorem]{Definition}

\newtheorem{eexample}[theorem]{Example}
\newenvironment{example}{\begin{eexample}\rm}{\end{eexample}}
\newtheorem{qquestion}[theorem]{Question}

\newtheorem{rremark}[theorem]{Remark}

\newcommand{\qed}{\rule{2.5mm}{2.5mm}}
\newenvironment{proof}{\begin{trivlist}\item[]{\sc Proof:}\quad}{%
         \hfill\qed\end{trivlist}}

\title{Pattern avoidance classes and subpermutations}
\author{M.D.~Atkinson\\
{\small{Department of Computer Science}}\\ {\small{University of Otago, New Zealand}}\\ 
\\M.M.~Murphy and N.~Ru\v skuc\\
{\small{School of Mathematics and Statistics}}\\
{\small{University of St Andrews, Scotland, KY16 9SS}}\\
}

\begin{document}
\maketitle
\begin{abstract} Pattern avoidance classes of permutations that cannot be expressed as unions of proper subclasses can be described as the set of subpermutations of a single bijection.  In the case that this bijection is a permutation of the natural numbers a structure theorem is given.  The structure theorem shows that the class is almost closed under direct sums or has a rational generating function.
\end{abstract}

\section{Introduction}\label{Intro}

Classes of permutations defined by their avoiding a given set of permutation 
patterns have been intensively studied within the last decade.  Quite often 
the issue has been to determine the number of permutations of each length 
in the class.  In order to do this it is necessary to derive structural 
properties of the permutations in the class starting from the avoided set.  
However, there are very few general techniques for obtaining such 
structural information.  This paper is a contribution towards a general 
structure theory.  We begin from the point of view that pattern-avoidance 
classes can be expressed as unions of {\em atomic} classes (those that have 
no non-trivial expression as a union).  We shall show that these atomic classes 
are precisely the classes that arise as the set of restrictions of 
some injection from one ordered set to another.  In general 
the order types of these two sets provide some information about the 
atomic class.  The major part of our paper is a characterisation of such 
injections and classes in the simplest case: when the order types are those 
of the natural numbers.

In the remainder of this section we review the terminology of pattern 
avoidance classes.  Most of this terminology is standard in the subject
except for the notion of a {\em natural} class.  We shall see  a number
of conditions on pattern avoidance classes that are equivalent to their being atomic and
we shall exhibit examples of atomic and non-atomic classes.  These
conditions and examples motivate the notion of a  
natural class whose elementary properties we explore in Section \ref{natural}.
Section \ref{main.theorem} contains our main result: a characterisation of 
natural classes, and Section \ref{example} gives some further examples.

We 
need a small number of definitions concerned with  permutations
and sets of permutations.  For our purposes a permutation is just an arrangement of the numbers 
$1,2,\ldots,n$ for some $n$.  We shall often need to consider arrangements 
of other sets of numbers and we shall refer to these as sequences; so, 
unless stated otherwise, a \emph{sequence} will mean
 a list of \emph{distinct} numbers.

Two finite sequences of the same length
 $\alpha = a_1a_2a_3\cdots$ and $\beta=b_1b_2b_3 \cdots$ are
said to be \emph{order isomorphic} (denoted as $\alpha\cong\beta$) if
for all $i,j$ we have $a_i<a_j$ if and only if $b_i<b_j$.
Any sequence defines a unique order isomorphic permutation; for example
$7496\cong 3142$.

A sequence $\alpha$ is said to be \emph{involved in} a 
sequence $\beta$ (denoted as $\alpha\inv\beta$) if
$\alpha$ is order isomorphic to a subsequence of $\beta$.
Usually, involvement is defined between permutations; for example 
$1324\inv 146325$ because of the subsequence $1435$.

It is easily seen that the involvement relation is a partial 
order on the set of all (finite) permutations.  We study it in 
terms of its order ideals which we call 
\emph{closed sets}.  A closed set $X$ of permutations has the defining
property that if
$\alpha\in X$ and $\delta\inv\alpha$ then $\delta\in X$.

Closed sets are most frequently specified by their \emph{basis}: the set
of permutations that are minimal subject to not lying in the closed set.
Once the basis $B$ is given the closed set is simply
\[\{\sigma\mid\beta\not\inv\sigma\mbox{ for all }\beta\in B\}\]
and we shall denote it by $\A(B)$.

Closed sets arise in the context of limited capability sorting
machines such as networks of stacks, queues and deques with a point of
input and a similar output.  Here the basis consists of minimal
sequences that cannot be sorted into some desirable order.  As
sequences can be sorted if and only if they do not involve any basis
elements, basis elements are frequently referred to as ``forbidden
patterns'' and the set of permutations that can be sorted by such
a mechanism is described as the set of all permutations that ``avoid''
that basis.

Much of the thrust of this paper is in specifying closed sets in a different
way.  Suppose that $A$ and $B$ are sets of real numbers and let $\pi$ be an 
injection from $A$ to $B$.  Then every finite subset 
$\{c_1,c_2,\ldots,c_n\}$ of $A$, where $c_1<c_2<\ldots <c_n$ maps
to a sequence $\pi(c_1)\,\pi(c_2)\,\ldots\,\pi(c_n)$ which is order
isomorphic to some permutation.  The set of permutations that arise in this
way is easily seen to be closed and we denote it by $\Sub(\pi:A\rightarrow B)$.
In many cases the domain and range of $\pi$ are evident from the context 
in which case we write simply $\Sub(\pi)$.  Also, since we may always replace
$B$ by the range of $\pi$ we shall, from now on, assume that $\pi$ is a bijection.

\begin{example}
Let $A=\{1-1/2^i,2-1/2^i\mid i=1,2,\ldots\}$ and
$B=\{1,2,\ldots\}$.  Let $\pi$ be defined by:

\[
\pi(x) = \left\{
\begin{array}{ll}
2i-1   & \mbox{if $x = 1-1/2^i$} \\
2i & \mbox{if $x = 2-1/2^i$}
\end{array}
    \right.
\]
Then it is easily seen that any finite increasing sequence of elements in $A$ maps
to an increasing sequence of odd integers followed by an increasing sequence of 
even integers.  From this it follows readily that the permutations of
$\Sub(\pi:A\rightarrow B)$ are precisely those that consist of two increasing segments.  As shown
in \cite{Restricted} this
closed set has basis $\{321, 3142, 2143\}$.  Notice that $A$ and $B$
have order
types $2\omega$ and $\omega$.  This particular closed set cannot be defined
as $\Sub(\pi:A\rightarrow B)$ with both $A$ and $B$ having order type
$\omega$. 
\end{example}

We now give several conditions on a closed set equivalent to it being expressible
as $\Sub(\pi:A\rightarrow B)$.

\begin{theorem}\label{atomic}
The following conditions on a closed set $X$ are equivalent:

\begin{enumerate}
\item $X=\Sub(\pi:A\rightarrow B)$ for some sets $A,B$ and bijection $\pi$.
\item $X$ cannot be expressed as a union of two proper closed subsets.
\item For every $\alpha,\beta\in X$ there exists $\gamma\in X$ such that
$\alpha\inv\gamma$ and $\beta\inv\gamma$.
\item $X$ contains permutations $\gamma_1\inv\gamma_2\inv\cdots$ such
that, for every $\alpha\in X$, we have $\alpha\inv\gamma_n$ for some $n$.
\end{enumerate}
\end{theorem}
\begin{proof}

$1\Rightarrow 2$.  Suppose $X=\Sub(\pi:A\rightarrow B)$ and yet there
exist proper closed subsets $Y,Z$ of $X$ such that $X=Y\cup Z$.  Then there
exist permutations $\rho\in X\setminus Y$ and $\sigma\in X\setminus Z$.
Therefore we can find subsequences $r_1<r_2<\cdots$ and $s_1<s_2<\cdots$ of
$A$ which are mapped by $\pi$ to subsequences order isomorphic to $\rho$ 
and $\sigma$.  The union of $\{r_1,r_2,\ldots\}$ with $\{s_1,s_2,\ldots\}$
defines a sequence $t_1<t_2<\cdots$ that is mapped by $\pi$ to a subsequence
order isomorphic to a permutation $\tau\in X$.  Obviously, $\rho\inv\tau$ and
$\sigma\inv\tau$.  However $\tau$ belongs to at least one of $Y$ or $Z$, say $\tau\in Y$.  Since $X$ is closed we have $\rho\in Y$, a contradiction.

$2\Rightarrow 3$.	Suppose that  there exist $\alpha,\beta\in X$ with the property 
that no permutation of $X$ involves both of them.  Put
	\begin{eqnarray*}
		Y&=&\{\gamma\in X\mid\alpha\not\inv\gamma\}\\
		Z&=&\{\gamma\in X\mid\beta\not\inv\gamma\} 
	\end{eqnarray*}
Then $Y$ and $Z$ are proper closed subsets of $X$ whose union is 
$X$ (since any $\gamma\in X\setminus (Y\cup Z)$ would involve both $\alpha$ 
and $\beta$).

$3\Rightarrow 4$.	If $\theta,\phi$ are two permutations in $X$ we know that
there exists a permutation of $X$ that involves both.  Temporarily we shall
use the notation $\theta\vee\phi$ to denote one of these permutations.  Now let $\beta_{1},\beta_{2},\ldots$ be any 
listing of the permutations of $X$.  We define a sequence of permutations of $X$
as follows: $\gamma_1=\beta_1$ and, for $i\geq 2$, $\gamma_i=\gamma_{i-1}
\vee\beta_i$.  Obviously, $\gamma_1\inv\gamma_2\inv\cdots$ and, for each permutation $\beta_n\in X$, $\beta_n\inv\gamma_n$.
	
$4\Rightarrow 1$.	In the sequence $\gamma_1\inv\gamma_2\inv\cdots$ we remove duplicates (if any) and we insert suitable permutations so that we have one of every degree.  This gives a 
sequence of permutations $\alpha_{1}\inv \alpha_{2}\inv\cdots$ 
such that
	\begin{enumerate}
		\item $|\alpha_{i}|=i$,
		\item $\alpha_{i}\in X$,
		\item for all $\sigma\in X$ there exists some $\alpha_{i}$ with 
		$\sigma\inv\alpha_{i}$
	\end{enumerate}
Now we shall inductively define, for each $i$, sets $A_{i},B_{i}$ 
and bijections $\pi_{i}:A_{i}\rightarrow B_{i}$ with the following 
properties:
	\begin{enumerate}
		\item $|A_{i}|=|B_{i}|=i$
		\item $\pi_{i}$ is order isomorphic to $\alpha_{i}$
		\item $A_{i-1}\subset A_{i}$ and $B_{i-1}\subset A_{i}$
		\item $\pi_{i}|_{A_{i-1}}=\pi_{i-1}$
	\end{enumerate}
Once these sets have been constructed we can complete the proof by 
setting $A=\bigcup_{i}A_{i}$ and $B=\bigcup_{i}B_{i}$.  Then we 
define $\pi:A\rightarrow B$ for any $a\in A$ by finding some $A_{i}$ for 
which $a\in A_{i}$ and setting $\pi(a)=\pi_{i}(a)$; by the last two 
properties $\pi$ is well-defined.  The second property guarantees 
that $X=\Sub(\pi:A\rightarrow B)$.

To carry out the construction we shall define $A_{i},B_{i}$ as 
subsets of the open interval $(0,1)$.  We
begin by setting
$A_{1}=B_{1}=1/2$ and $\pi_1(1/2)=1/2$.  Suppose now that 
$A_{i},B_{i}, \pi_{i}$ have been constructed for $i=1,2\ldots,n$.  
The permutation $\alpha_{n+1}$ is constructed from $\alpha_{n}$ by 
the insertion of a new element $t$ at position $s$ in $\alpha_{n}$; 
the position numbers of all the elements of $\alpha_{n}$ which are 
greater than or equal to $s$ have to be increased by $1$ and those 
values which are greater than or equal to $t$ have also to be 
incremented by $1$.

We reflect this insertion in the definition of $A_{n+1},B_{n+1}$ and 
$\pi_{n+1}$.  The set $A_{n+1}$ is formed by augmenting $A_{n}$ with 
another number $a$ whose value lies between its $(s-1)^{th}$ and $s^{th}$ elements
(if $s=1$ we take $a$ between $0$ and the minimal element of $A_n$; while if
$s=n$ we take $a$ between the maximal element and $1$).  Similarly 
 $B_{n+1}$ 
is formed by augmenting $B_{n}$ with 
a number $b$ whose value lies between its $(t-1)^{th}$ and $t^{th}$ elements.  Then we 
define $\pi_{n+1}$ so that it agrees with $\pi_{n}$ on the elements 
of $A_{n}$ and has $\pi_{n+1}(a)=b$.
\end{proof}

Because of this result we call closed sets of the form 
$\Sub(\pi:A\rightarrow B)$ \emph{atomic} on the grounds that
they cannot be decomposed as a proper union of two closed
subsets.  Expressing a given closed set as a union of atomic sets
is often very useful in discovering structural information.

\begin{example}(See \cite{Restricted})
$\A(321,2143)=\A(321,2143,3142)\cup \A(321,2143,2413)$
\end{example}

Given an arbitrary closed subset one might hope
to find its properties by first expressing it as a union of
atomic sets, and then discovering properties of the bijection
$\pi$ associated with each atomic subset.  Many difficulties
impede this approach.  It may happen that a closed set cannot
be expressed as a finite union of atomic subsets.  Moreover
an atomic closed set may have a defining bijection $\pi$
whose domain and range have high ordinal type; in that case one might
be hopeful that properties of these ordinals (in particular, limit points)
might imply properties of $X=\Sub(\pi:A\rightarrow B)$.  Despite
this hope it seems sensible
to begin the systematic study of atomic sets by looking at 
the case where the ordinal type of both $A$ and
$B$ is that of the natural numbers $\Nat$.

\section{Natural classes and sum-complete classes}\label{natural}

A \emph{natural class} is a closed set of the form
$\Sub(\pi:\Nat\rightarrow \Nat)$.  In other words, starting from
a permutation  $\pi$
of the natural numbers, we form all the finite subsequences of
$\pi(1),\pi(2),\ldots$ and define a natural class as consisting of the permutations
order isomorphic to these subsequences.  From now on we shall use
the notation $\Sub(\pi)$ (suppressing a notational reference to the
domain and range of $\pi$) in the following circumstances
\begin{enumerate}
\item when $\pi$ is an infinite permutation with $\Nat$ as its domain and range,
\item when $\pi$ is a finite permutation (in which case the domain and range
are $\{1,2,\ldots,n\}$ where $n$ is the degree of $\pi$).
\end{enumerate}

\begin{example}
Let $\pi$ be defined by:
\[
\pi = 1\; 3\; 2\; 6\; 5\; 4\; 10\; 9 \; 8\; 7 \; \ldots
\]
Then $\Sub(\pi)$ is easily seen to be the set of all layered permutations as defined
in \cite{layered}.
\end{example}

If $\alpha = a_1a_2\cdots a_m$ and
$\beta=b_1b_2\cdots $ are sequences (in particular, permutations)
then their \emph{sum} $\alpha\oplus
\beta$ is defined to be the permutation $\gamma\delta$ where
the segments $\gamma$ and $\delta$ are rearrangements of
$1,2,\ldots,m$ and $m+1,m+2,\ldots$ respectively, and 
$\alpha\cong\gamma$ and $\beta\cong\delta$.  Notice that we
do not require that $\beta$ be a finite permutation.
If a permutation can
be expressed as $\alpha\oplus\beta$ (with neither summand empty) we
say that it is \emph{decomposable}; otherwise 
it is said to be \emph{indecomposable}.  We also extend the sum
notation to sets by defining, for any two sets of permutations $X$ and $Y$,
\[X\oplus Y=\{\sigma\oplus\tau\mid \sigma\in X, \tau\in Y\}\]

A set $X$ of permutations is said to be \emph{sum-complete}
if for all $\alpha, \beta\in X$, we have
$\alpha\oplus\beta\in X$.

Sum-completeness and decomposability are linked by
the
following result, proved in \cite{Wreath}.

\begin{lemma}\label{sum.comp.indecomp.lemma}
Let $X$ be a closed set with basis $B$.  Then $X$ is sum-complete if 
and only if $B$ contains only indecomposable permutations.
\end{lemma}

We shall see that natural classes and sum-complete closed sets
are closely connected.  The first hint of this connection is the following result
which, in particular, shows that
every sum-complete closed set is a natural class.

\begin{proposition}\label{most.natural}
Let $\gamma$ be any (finite) permutation and $S$
any sum-complete closed set.  Then $\Sub(\gamma)\oplus S$ is a natural
class.
\end{proposition}
\begin{proof}Let $\beta_1,\beta_2,\ldots$ be any listing of the permutations
of $S$.  Consider the sequence of permutations
\[\gamma\inv\gamma\oplus\beta_1\inv\gamma\oplus\beta_1\oplus\beta_2
\inv\gamma\oplus\beta_1\oplus\beta_2\oplus\beta_3\inv\cdots\]
Since $S$ is sum-complete all these permutations lie in 
$\Sub(\gamma)\oplus S$.  On the other hand it is clear that every permutation
of $\Sub(\gamma)\oplus S$ is involved in some term of the sequence.  Hence
 $\Sub(\gamma)\oplus S$ satisfies
condition 4 of Theorem \ref{atomic}, and hence is atomic.
Furthermore, the proof of (4$\Rightarrow$1) in Theorem \ref{atomic} tells us how to express
$\Sub(\gamma)\oplus S$ in the form $\Sub(\pi:A\rightarrow B)$.
Following this recipe, it is easy to see that both $A$ and $B$ are (isomorphic to) $\mathbb{N}$,
and we have
a natural class, as required.
\end{proof}

Notice that the proof of this result makes no assumption on the listing
of the elements of $S$.  That means that the infinite permutation $\pi$ for which
$\Sub(\gamma)\oplus S=\Sub(\pi)$ is very far from being unique.

In the remainder of the paper we shall be exploring a partial converse of
Proposition \ref{most.natural}.  Our main theorem will show
that every finitely based natural class $X$ does have the form of the proposition
unless $\pi$ and $X$ have a very particular form.

\section{A characterisation of natural classes}\label{main.theorem}

This section is devoted to the proof of the following theorem.

\begin{theorem}\label{main}
Let $X$ be a finitely based natural class.  Then either
\begin{enumerate}
\item $X=\Sub(\gamma)\oplus S$ where $\gamma$ is a finite permutation
and $S$ is a sum-complete closed class determined uniquely by $X$, or
\item $X=\Sub(\pi)$ where $\pi$ is unique and ultimately periodic
in the sense that there exist integers $N$ and $P>0$ such that, for all $n\geq N$,
$\pi(n+P)=\pi(n)+P$.
\end{enumerate}
\end{theorem}

The proof of the theorem will show precisely how $X$ determines $S$ in the
first alternative.  It will also, in the case of the second alternative, prove that
$X$ is enumerated by a rational
generating function.

Before embarking on a series of lemmas that lead up to the proof of 
Theorem \ref{main} we shall define some notation that will be in force for the rest
of this section.

We shall let $X=\Sub(\pi)$ where $\pi$ is a permutation of $\Nat$.  The basis
of $X$ will be denoted by $B$ and we let $b$ denote the length of a longest
permutation in $B$.  The permutations of $B$ have a decomposition into
sum components; the set of final components in such decompositions will be 
denoted by $C$.

We shall use the notation $\A(C)$ for the closed set of all permutations that
avoid the permutations of $C$.  This is a slight extension of the notation we defined
in Section \ref{Intro} because $C$ might not be the basis of $\A(C)$ ($C$ might
contain some non-minimal elements outside $\A(C)$).  This causes no
technical difficulties.  Obviously, as every permutation that avoids the permutations of
$C$ also avoids the permutations of $B$, we have $\A(C)\subseteq X$.  By Lemma \ref{sum.comp.indecomp.lemma} 
$\A(C)$ is sum-complete; it is, as we shall see,
the sum-complete class $S$ occurring in the statement
of Theorem \ref{main}.

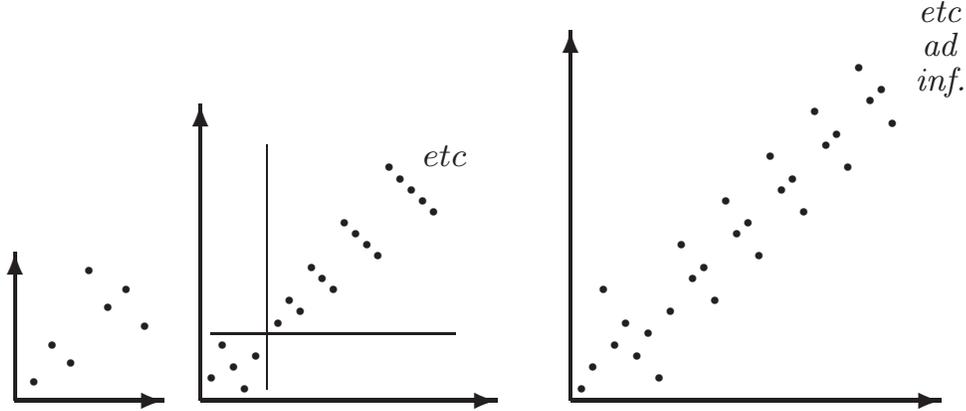
\begin{figure}
\begin{center}
\scalebox{1.4}{ 
\begin{picture}(250,110)(0,0)
\put(0,0){\makebox{  
\thicklines
\put(0,0){\vector(1,0){40}}
\put(0,0){\vector(0,1){40}}
\put(5,5){\circle*{2}}
\put(10,15){\circle*{2}}
\put(15,10){\circle*{2}}
\put(20,35){\circle*{2}}
\put(25,25){\circle*{2}}
\put(30,30){\circle*{2}}
\put(35,20){\circle*{2}}

}}

\put(50,0){\makebox{  
\thicklines
\put(0,0){\vector(1,0){80}}
\put(0,0){\vector(0,1){80}}
\put(3,6){\circle*{2}}
\put(6,15){\circle*{2}}
\put(9,9){\circle*{2}}
\put(12,3){\circle*{2}}
\put(15,12){\circle*{2}}

\thinlines
\put(3,18){\line(1,0){66}}
\put(18,3){\line(0,1){66}}

\put(21,21){\circle*{2}}

\put(24,27){\circle*{2}}
\put(27,24){\circle*{2}}

\put(30,36){\circle*{2}}
\put(33,33){\circle*{2}}
\put(36,30){\circle*{2}}
\put(39,48){\circle*{2}}
\put(42,45){\circle*{2}}
\put(45,42){\circle*{2}}
\put(48,39){\circle*{2}}
\put(51,63){\circle*{2}}
\put(54,60){\circle*{2}}
\put(57,57){\circle*{2}}
\put(60,54){\circle*{2}}
\put(63,51){\circle*{2}}
\put(66,66){\makebox(0,0){\small{$etc$}}}
}}

\put(150,0){\makebox{  
\thicklines
\put(0,0){\vector(1,0){100}}
\put(0,0){\vector(0,1){100}}
\put(3,3){\circle*{2}}
\put(6,9){\circle*{2}}
\put(9,30){\circle*{2}}
\put(12,15){\circle*{2}}
\put(15,21){\circle*{2}}
\put(18,12){\circle*{2}}
\put(21,18){\circle*{2}}
\put(24,6){\circle*{2}}
\put(27,24){\circle*{2}}

\put(30,42){\circle*{2}}
\put(33,33){\circle*{2}}
\put(36,36){\circle*{2}}
\put(39,27){\circle*{2}}

\put(42,54){\circle*{2}}
\put(45,45){\circle*{2}}
\put(48,48){\circle*{2}}
\put(51,39){\circle*{2}}

\put(54,66){\circle*{2}}
\put(57,57){\circle*{2}}
\put(60,60){\circle*{2}}
\put(63,51){\circle*{2}}

\put(66,78){\circle*{2}}
\put(69,69){\circle*{2}}
\put(72,72){\circle*{2}}
\put(75,63){\circle*{2}}

\put(78,90){\circle*{2}}
\put(81,81){\circle*{2}}
\put(84,84){\circle*{2}}
\put(87,75){\circle*{2}}

\put(100,95){\makebox(0,0){\small{\shortstack{\emph{etc}\\ \emph{ad}\\ 
\emph{inf.}}}}}

}}

\end{picture}
}  
\caption{The sum of $132$ and $4231$ is $132\,7564$, as plotted on 
the left.  Every finitely based natural class is defined by a 
finite permutation summed with a sum-complete class (centre), or is 
eventually periodic (right). \label{sum.fig}}
\end{center}
\end{figure}

From time to time we shall illustrate our proof with 
diagrams that display permutations.
These diagrams are plots in the $(x,y)$ plane.  A permutation $p_1,p_2,\ldots$
(which maps $i$ to $p_i$) will be represented by a set of points whose
coordinates are $(i,p_i)$.  As a first use of such diagrams we have 
Figure \ref{sum.fig}
which illustrates the sum operation and the two alternatives in 
Theorem \ref{main}.

\begin{lemma}\label{k-exists}There exists an integer $k$ such that, for all $d>k$,
\[
\Sub(\pi(d),\pi(d+1),\ldots)=\A(C)
\]
\end{lemma}
\begin{proof}
For each $\gamma\in C$ there is a basis element of $X$ of the form 
$\beta\oplus\gamma$.  Every such $\beta$ is a permutation of $X$
and so we can choose a subsequence $S(\beta)$ of $\pi$ with 
$S(\beta)\cong\beta$.  Let $t$ be the maximal value occurring in 
all such $S(\beta)$ and let $u$ be the right-most position of $\pi$
where an element of some $S(\beta)$ occurs.

There exists an integer $k>u$ such that all terms $\pi(k+1),\pi(k+2),\ldots$
exceed $t$.  Note that the order type of $\Nat$ is used in establishing the
existence of $k$.  Among the terms $\pi(k+1),\pi(k+2),\ldots$
there can be no subsequence order isomorphic to an element of $C$.
This proves that $\Sub(\pi(d),\pi(d+1),\ldots)\subseteq\A(C)$ for
all $d>k$.  It also proves that $\A(C)$ is non-empty.

Now let $\theta\in\A(C)$.  Since the permutation $1$ lies in $\A(C)$
and $\A(C)$ is sum-complete we have $1,2,\ldots, d-1\oplus \theta\in\A(C)$.
Therefore $\pi$ has a subsequence order isomorphic to this permutation
and that implies that $\pi(d),\pi(d+1),\ldots$ has a subsequence order
isomorphic to $\theta$ which completes the proof.
\end{proof}

\begin{corollary}\label{cases}
Either $X=\Sub(\gamma)\oplus \A(C)$ for some finite permutation $\gamma$,
or $\pi$ has finitely many components and the last component (which is necessarily infinite)
involves an
element of $C$.
\end{corollary}
\begin{proof}
Let $\pi=\pi_1\oplus\pi_2\oplus\cdots$ be the sum decomposition of $\pi$.
Lemma \ref{k-exists}
tells us, in particular, that there is a maximal position $k$ where
a subsequence order isomorphic to an element of $C$ can begin.  Suppose this
position occurs in the sum component $\pi_r$.  If $\pi_r$ is not the final component
of $\pi$ then we have 
$\Sub(\pi)=\Sub(\pi_1\oplus\cdots\oplus\pi_r)\oplus\Sub(\pi_{r+1}\oplus\cdots)$.
However $\gamma=\pi_1\oplus\cdots\oplus\pi_r$ is finite and 
$\Sub(\pi_{r+1}\oplus\cdots)=\A(C)$ by the lemma.
\end{proof}

The first alternative of this corollary leads to the first alternative of Theorem 
\ref{main} because of the following uniqueness result.

\begin{proposition}
If $X=\Sub(\gamma_1)\oplus S_1=\Sub(\gamma_2)\oplus S_2$ where 
$\gamma_1,\gamma_2$ are finite permutations and $S_1,S_2$ are sum-complete
then $S_1=S_2$.
\end{proposition}
\begin{proof}
Let $\sigma_1\in S_1$.  Then, as $S_1$ contains every permutation of the form
$\iota_m=1\ 2\ \ldots m$, $S_1$ also contains $\iota_t\oplus\sigma_1$ where
$t=|\gamma_2|$.  But this permutation belongs to $\Sub(\gamma_2)\oplus S_2$ and
so can be expressed as $\gamma'\oplus\sigma_2$ where $\gamma'\inv\gamma_2$.
Since $\iota_t\oplus\sigma_1=\gamma'\oplus\sigma_2$ and $|\gamma'|\leq|\iota_t|$
we have $\sigma_1\inv\sigma_2$.  This proves that 
$\sigma_1\in S_2$ and therefore $S_1\subseteq S_2$. The
result now follows by symmetry.
\end{proof}



In the remainder of the proof of Theorem \ref{main} we shall assume that
the second alternative of Corollary \ref{cases} holds and work towards
proving the second alternative of the theorem.  In particular, there
exists a greatest position $k$ in $\pi$ where a subsequence isomorphic
to a permutation in $C$ can begin, and this position occurs in the final (infinite)
sum component $\pi_z$ of $\pi$.

Next we prepare the ground for two arguments that occur later in the proof
and which depend upon the indecomposability of $\pi_z$.  Suppose that
$r$ is any position in $\pi_z$.  We define a pair of sequences 
$U(r)=u_1u_2\cdots$ and $V(r)=v_1v_2\cdots$ by the following rules:
\begin{enumerate}
\item $v_i$ is the position among the
terms of $\pi$ up to and including position $u_{i-1}$ (when $i=1$ take $u_0=r$)
where the greatest 
element occurs:
$$
\pi(v_i)=\max\{\pi(k) \:|\: k\leq u_{i-1}\}.
$$
\item 
$u_i$ is the rightmost position in $\pi$ where a term not exceeding
$\pi(v_i)$ occurs:
$$
u_i=\max\{k\:|\: \pi(k)\leq \pi(v_i)\}.
$$
\end{enumerate}

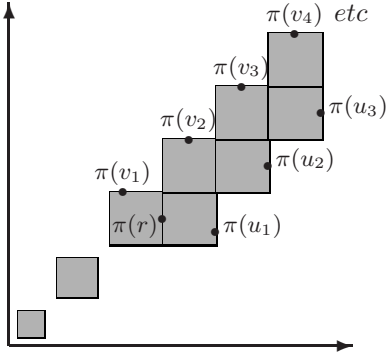
\begin{figure}
\begin{center}
\begin{picture}(150,120)(0,0)
\definecolor{rusty}{gray}{0.7}
\thicklines
\put(2,2){\vector(1,0){130}}
\put(2,2){\vector(0,1){130}}
\thinlines
\newlength{\origboxborder}
\setlength{\origboxborder}{\fboxsep}
\setlength{\fboxsep}{0mm}

\put(5,5){\fcolorbox{black}{rusty}{\makebox(10,10){}}}
\put(20,20){\fcolorbox{black}{rusty}{\makebox(15,15){}}}

\put(40,40){\fcolorbox{black}{rusty}{\makebox(20,20){}}}

\put(60,40){\fcolorbox{black}{rusty}{\makebox(20,20){}}}
\put(60,60){\fcolorbox{black}{rusty}{\makebox(20,20){}}}

\put(80,60){\fcolorbox{black}{rusty}{\makebox(20,20){}}}
\put(80,80){\fcolorbox{black}{rusty}{\makebox(20,20){}}}

\put(100,80){\fcolorbox{black}{rusty}{\makebox(20,20){}}}
\put(100,100){\fcolorbox{black}{rusty}{\makebox(20,20){}}}

\put(45,60){\circle*{3}}
\put(60,50){\circle*{3}}
\put(45, 65){\makebox[0pt][c]{\small{$\pi(v_1)$}}}
\put(41, 45){\makebox[0pt][l]{\small{$\pi(r)$}}}

\put(70,80){\circle*{3}}
\put(80,45){\circle*{3}}
\put(70, 85){\makebox[0pt][c]{\small{$\pi(v_2)$}}}
\put(83, 45){\makebox[0pt][l]{\small{$\pi(u_1)$}}}

\put(90,100){\circle*{3}}
\put(100,70){\circle*{3}}
\put(90, 105){\makebox[0pt][c]{\small{$\pi(v_3)$}}}
\put(103, 70){\makebox[0pt][l]{\small{$\pi(u_2)$}}}

\put(110,120){\circle*{3}}
\put(120,90){\circle*{3}}
\put(110,125){\makebox[0pt][c]{\small{$\pi(v_4)$}}}
\put(123,90){\makebox[0pt][l]{\small{$\pi(u_3)$}}}

\put(125,125){\makebox{$etc$}}
\setlength{\fboxsep}{\origboxborder} 
\end{picture}  %
\end{center}
\caption{The terms of $\pi$ as mapped out by $\pi(u_i)$ and 
$\pi(v_i)$.  All terms lie in the shaded boxes.}\label{pi.fig}
\end{figure}

Figure \ref{pi.fig} depicts these points and the next lemma assures us that
the figure accurately represents the relative positions of the marked points.

\begin{lemma}
The relative positions and sizes of the terms $\pi(u_i)$ and 
$\pi(v_j)$ are described by the following inequalities:
\[
v_1 < v_2 < u_1 < v_3 < u_2 < v_4 < u_3 < \cdots
\]
$$
\pi(u_1)<\pi(v_1) < \pi(u_2) < \pi(v_2) < \pi(u_3) < 
\pi(v_3) < \cdots
$$
\end{lemma}

\begin{proof}
(I) From the definition of $v_i$ we have $v_i\leq u_{i-1}$, and from the definition
of $u_i$ we have 
and $\pi(v_i)\geq \pi(u_{i})$.
Note that we cannot have $u_i=u_{i-1}$, because then every term of $\pi$ to the left of this position
would be less than or equal to $\pi(v_i)$,
and every term to the right would be greater than $\pi(v_i)$,
contradicting the assumption that $\pi(v_i)$ belongs to the final component of $\pi$.
Hence we have $v_i\leq u_{i-1}<u_i$ and $\pi(v_i)>\pi(u_i)$.

(II)
By the definition of $v_{i+1}$, we have $v_{i+1}\leq u_i$ and $\pi(v_{i+1})\geq \pi(v_i)$.
We cannot have $v_{i+1}< v_i$, because $\pi(v_i)$ is the maximal value of $\pi$
on the interval $[1,u_{i-1}]$.
Also, we cannot have $v_{i+1}=v_i$, because that would imply $u_{i+1}=u_i$,
which is proved impossible as in (I).
Finally, we cannot have $v_{i+1}=u_i$ because $\pi(u_i)<\pi(v_i)$ by (I).
We conclude that $v_i<v_{i+1}<u_i$ and $\pi(v_i)<\pi(v_{i+1})$.

(III)
As in (I), we have $u_{i+1}>u_i$ and $\pi(u_{i+1})<\pi(v_{i+1})$.
Moreover, $u_{i+1}>u_i$ immediately implies that $\pi(u_{i+1})>\pi(v_i)$.

(IV)
As in (II), we have $v_{i+2}<u_{i+1}$ and $\pi(v_{i+2})>\pi(v_{i+1})$.
Moreover, $v_{i+2}>u_i$ for otherwise we would have
$\pi(v_{i+2})\leq \pi(v_{i+1})$.

Summarising (I)--(IV), we have
$v_i<v_{i+1}<u_i<v_{i+2}<u_{i+1}$ and
$\pi(u_i)<\pi(v_i)<\pi(u_{i+1})<\pi(v_{i+1})<\pi(v_{i+2})$ for every $i=1,2,\ldots$,
which is enough to prove the lemma.
\end{proof}

Our first use of the sequences $U(r)$ and $V(r)$ and 
the above lemma occurs immediately.  We have
seen (Lemma \ref{k-exists}) that there is a rightmost position in $\pi$ where
subsequences order isomorphic to permutations in $C$ can begin.  Now
we prove that there is a rightmost position by which they have all ended.

\begin{lemma}\label{l-exists}There exists a position $\ell$ of $\pi$ such that no
subsequence of $\pi$ that is order isomorphic to an element of $C$
terminates after position $\ell$.
\end{lemma}
\begin{proof}Consider the sequences $U(k),V(k)$ and refer to Figure 
\ref{pi.fig} with $r=k$, 
in particular to the edge-connected strip of boxes that
begins with the box $B_1$ bounded by $\pi(v_1)$ and $\pi(k)$. 
Let $S$ be a subsequence of $\pi$ isomorphic to
a permutation $\gamma\in C$.  
By definition of $k$, $S$ cannot start to the right of $\pi(k)$.
In fact, since $\gamma$ is indecomposable, $S$ must start in $B_1$, and the
terms of $S$ must lie in a contiguous  segment of boxes. 
 Therefore,
as $|S|\leq b$, $S$ cannot extend beyond position $u_{b/2}$.
\end{proof}

In view of this lemma we may define $\ell$ as the last position of $\pi$ that
is part of a subsequence isomorphic to an element of $C$.  Now we define
the sequences $U(\ell), V(\ell)$ (and, re-using notation, call them
$u_1,u_2,\ldots$ and $v_1,v_2,\ldots$).

The defining property of $\ell$ implies that $\pi(1),\ldots,\pi(\ell)$
is the only subsequence of $\pi$ with this order isomorphism type.  For
any subsequence of $\pi$ order isomorphic to 
$\pi(1),\ldots,\pi(\ell)$ has a final element that is part of a subsequence
order isomorphic to an element of $C$.  Therefore this final element cannot occur
after position $\ell$ within $\pi$ and so it must be $\pi(1),\ldots,\pi(\ell)$
itself.  Notice also, again from the definition of $\ell$, that the permutation
order isomorphic to $\pi(1),\ldots,\pi(\ell)$ is the longest permutation
in $X$ whose last element is the terminating element of a subsequence
order isomorphic to an element of $C$; as such, this permutation
depends on $X$ and not on $\pi$.
In the next lemma we prove that a number of other initial segments
of $\pi$ are unique of their isomorphism type, and depend on $X$ rather than on
$\pi$.

\begin{lemma}\label{uniqueness1}
For each $i=1,2,\ldots$ the sequence $\pi(1),\ldots,\pi(u_i)$ is  the unique
subsequence of that order isomorphism type.  
Its corresponding permutation $\beta=\beta(1),\ldots,\beta(u_i)$ is
the longest permutation in $X$ satisfying the following two properties:
\begin{itemize}
\item[{\rm (1)}]
$\beta(1),\ldots,\beta(u_{i-1})$ is isomorphic to $\pi(1),\ldots,\pi(u_{i-1})$
(where $u_0=l$);
\item[{\rm (2)}]
$\beta(u_i)\leq \beta(v_i)$.
\end{itemize}
As such, $\beta$ depends on $X$ only, and not on $\pi$.
\end{lemma}

\begin{proof}
We prove the lemma by induction, anchoring it at $u_0=l$.
Assume that the statements are true for some $i\geq 1$, and consider
any subsequence $\pi(s_1),\ldots,\pi(s_{u_{i+1}})$ order isomorphic to
$\pi(1),\ldots,\pi(u_{i+1})$.
By the inductive hypothesis we must have $s_j=j$ for $j=1,\ldots,u_i$.
But then, since $\pi(u_{i+1})$ is the rightmost term of $\pi$
smaller than $\pi(v_{i+1})$,
and since there are $u_{i+1}-u_i-1$  terms between $\pi(u_i)$ and $\pi(u_{i+1})$,
it follows that $s_j=j$ for $j=u_i+1,\ldots,u_{i+1}$ as well.

Clearly, the permutation $\beta$ satisfies properties (1) and (2), by virtue of being isomorphic to
$\pi(1),\ldots,\pi(u_{i+1})$.
Suppose that $\gamma=\gamma(1),\ldots,\gamma(m)$ is any permutation satisfying these conditions.
Consider an embedding $\pi(t_1),\ldots,\pi(t_m)$ of $\gamma$ in $\pi$.
As above, we must have $t_j=j$ for $j=1,\ldots,u_i$.
And again, $\pi(u_{i+1})$ being the rightmost term of $\pi$ smaller than $\pi(v_{i+1})$,
we have that $t_m\leq u_{i+1}$.
But this, in turn, implies that $|\gamma|=m\leq t_m\leq u_{i+1}=|\beta|$.
This proves that $\beta$ is indeed the longest permutation of $X$ satisfying (1) and (2).
The last statement of the lemma is now straightforward.
\end{proof}

At this point we can prove the uniqueness of $\pi$:  it is the limit
of its initial segments $\pi(1)\cdots\pi(u_i)$ and these depend on $X$
alone.

For future use we record the following result, the proof of which is analogous to the proof of
Lemma \ref{uniqueness1}:

\begin{lemma}\label{uniqueness2}
For each $i=2,3,\ldots$ the subsequence consisting of all terms of $\pi$ not exceeding
$\pi(v_i)$ is unique of its order isomorphism type.  
Its corresponding permutation $\beta=\beta(1),\ldots,\beta(n_i)$ is
the longest permutation in $X$ satisfying the following two properties:
\begin{itemize}
\item[{\rm (1)}]
$\beta(1),\ldots,\beta(u_{i-1})$ is isomorphic to $\pi(1),\ldots,\pi(u_{i-1})$;
\item[{\rm (2)}]
$\beta(v_i)$ is its largest term.
\end{itemize}
As such, $\beta$ depends on $X$ only, and not on $\pi$.
\end{lemma}

The subsequences whose embeddings are unique in the previous two lemmas
are those all of whose terms are taken from an initial contiguous strip
of blocks in Figure \ref{pi.fig} (with $r=\ell$).  For convenience we
let $\sigma(i)$ be the permutation isomorphic to the first type
of subsequence (Lemma \ref{uniqueness1}), and $\sigma'(i)$ be the permutation isomorphic to the
second type (Lemma \ref{uniqueness2}).

\begin{figure}[b!]
\begin{center}
\begin{picture}(150,120)(0,0)
\definecolor{rusty}{gray}{0.7}
\thicklines
\put(2,2){\vector(1,0){130}}
\put(2,2){\vector(0,1){130}}
\thinlines
\setlength{\origboxborder}{\fboxsep}
\setlength{\fboxsep}{0mm}

\put(5,5){\fcolorbox{black}{rusty}{\makebox(30,30){}}}

\put(40,40){\fcolorbox{black}{rusty}{\makebox(20,20){}}}

\put(60,40){\fcolorbox{black}{rusty}{\makebox(20,20){}}}
\put(60,60){\fcolorbox{black}{rusty}{\makebox(20,20){}}}

\put(80,60){\fcolorbox{black}{rusty}{\makebox(20,20){}}}
\put(80,80){\fcolorbox{black}{rusty}{\makebox(20,20){}}}

\put(100,80){\fcolorbox{black}{rusty}{\makebox(20,20){}}}
\put(100,100){\fcolorbox{black}{rusty}{\makebox(20,20){}}}

\put(45,60){\circle*{3}}
\put(60,50){\circle*{3}}
\put(45, 65){\makebox[0pt][c]{\small{$\pi(v_1)$}}}
\put(41, 45){\makebox[0pt][l]{\small{$\pi(l)$}}}

\put(70,80){\circle*{3}}
\put(80,45){\circle*{3}}
\put(70, 85){\makebox[0pt][c]{\small{$\pi(v_2)$}}}
\put(83, 45){\makebox[0pt][l]{\small{$\pi(u_1)$}}}

\put(90,100){\circle*{3}}
\put(100,70){\circle*{3}}
\put(90, 105){\makebox[0pt][c]{\small{$\pi(v_3)$}}}
\put(103, 70){\makebox[0pt][l]{\small{$\pi(u_2)$}}}

\put(110,120){\circle*{3}}
\put(120,90){\circle*{3}}
\put(110,125){\makebox[0pt][c]{\small{$\pi(v_4)$}}}
\put(123,90){\makebox[0pt][l]{\small{$\pi(u_3)$}}}

\put(65,45){\makebox{$L$}}
\put(65,65){\makebox{$U$}}
\put(85,65){\makebox{$L$}}
\put(85,85){\makebox{$U$}}
\put(105,85){\makebox{$L$}}
\put(105,105){\makebox{$U$}}

\put(125,125){\makebox{$etc$}}
\setlength{\fboxsep}{\origboxborder} 
\end{picture}  %
\end{center}
\caption{$\pi$ mapped out by $u_i$ and $v_i$.  None of the boxes marked 
$L$ or $U$ contains a monotonic increasing or decreasing subsequence 
of length $b$.  Thus the sizes of these boxes is bounded.}\label{uvpi.fig}
\end{figure}
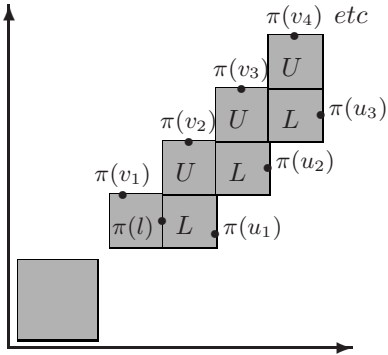

\begin{lemma}\label{bounded.separation}$u_{i+1}-u_i\leq 2(b-1)^2$
\end{lemma}
\begin{proof}
Of course $u_{i+1}-u_i$ is the number of terms of 
$\pi(u_{i}+1)\cdots\pi(u_{i+1})$.  We divide these terms into two sets
\[L=\{j\mid u_i<j\leq u_{i+1}\mbox { and }\pi(j)<\pi(v_{i+1})\}\]
and
\[U=\{j\mid u_i<j\leq u_{i+1}\mbox { and } \pi(j)>\pi(v_{i+1})\}\]
and we shall show that both $|L|$ and $|U|$ are at most $(b-1)^2$.
Each bound is proved in the same way and we give the details for $|L|$ only.
Figure \ref{uvpi.fig} depicts the locations of $L$ and $U$ within $\pi$.

Consider a maximal increasing subsequence $\pi(j_1)\pi(j_2)\cdots\pi(j_m)$
of $\pi$ such that all $j_i\in L$.   Using this we form another subsequence
of $\pi$ whose terms are the following:
\begin{enumerate}
\item all terms not exceeding $\pi(v_i)$, (a subsequence order isomorphic to
$\sigma'(i)$)
\item the term $\pi(v_{i+1})$
\item all the terms $\pi(j_1)\pi(j_2)\cdots\pi(j_m)$
\end{enumerate}

The permutation which is order isomorphic to this subsequence is, of course,
a member of $X$ and we write it as $\lambda n \rho\theta$ where $n$ corresponds
to the term $\pi(v_{i+1})$, $\lambda\rho$ corresponds to $\sigma'(i)$, and 
$\theta=n-m,n-m+1,\ldots,n-1$ corresponds to $\pi(j_1)\pi(j_2)\cdots\pi(j_m)$.

Now consider another permutation almost the same as this except that $\theta$ 
contains one more term, and $n$ is replaced by $n+1$.  It has the form
$\lambda, n+1, \rho\theta'$ with $\theta'=n-m,n-m+1,\ldots,n-1,n$.  This permutation
does not belong to $X$.  To see this, assume that some subsequence of $\pi$
is order isomorphic to it.  In the correspondence between the permutation and the subsequence,
$\lambda\rho$ (which is order isomorphic to $\sigma'(i)$) must be mapped
to the subsequence of terms not exceeding $\pi(v_i)$, 
by Lemma \ref{uniqueness2}, and $n+1$ must be mapped to one
of the terms of $\pi$ in the range of positions $u_{i-1}+1$ to $u_i$.  This forces
$\theta'$ to be mapped into $L$ as these are the only positions of $\pi$ to the right
of $u_i$ and smaller than $\pi(v_{i+1})$.  However, this contradicts that $L$
contains no increasing sequence of length $m+1$.

It follows that $\lambda, n+1, \rho\theta'$ must involve a basis element of $X$.
A particular embedding of a basis element must contain \emph{all} the terms of
$\theta'$ for otherwise this basis element would be embedded in $\lambda n \rho\theta$ which is impossible.  In particular we can deduce that $m+1\leq b$.

Exactly the same argument can be carried out for maximal decreasing
subsequences.  Thus the sequence $\pi(L)$ contains no
increasing or decreasing subsequence of length more than $b-1$ and, by the
well known result of Erd\H os and Szekeres (see \cite{Cameron}), we conclude
that $|L|\leq (b-1)^2$.

The proof that $|U|\leq (b-1)^2$ is similar but it uses 
$\pi(u_{i+1})$ and  $\sigma(i)$ instead of $\pi(v_{i+1})$ 
and  $\sigma'(i)$.
\end{proof}

We now define an encoding of permutations of $X$.  
If $\xi=x_1x_2\cdots x_n\in X$ then we encode it as 
$E(\xi)=e_1e_2\cdots e_n$ where
\[e_k=|\{i\mid 1\leq i\leq k \mbox{ and }x_i\geq x_k\}| \]

If a term of $\pi$ lies in the final component beyond position $u_2$
then there are, by
Lemma \ref{bounded.separation}, at most $4(b-1)^2$ preceding terms
greater than it.  On the other hand, if it lies in one of the finite components or
in the final component and not beyond position $u_2$ then, obviously,
there will again only be a bounded number of preceding greater terms.  Since
permutations of $X$ are order isomorphic to subsequences of $\pi$
and $\pi$ is determined by $X$ there is
an upper bound depending on $X$ alone for each of the code symbols $e_k$.
Thus we may consider $E(X)$ as a language over some finite alphabet $A=\{1,\ldots,m\}$.

The encoding has the property that, if $\gamma=\delta\epsilon$ is a permutation
and $\bar{\delta}$ is the permutation order isomorphic to the initial segment
$\delta$, then $E(\bar{\delta})$ is an initial segment of $E(\gamma)$.

Furthermore this encoding  (in which every element of a permutation is encoded
by the number of its higher predecessors) is closely related to the encoding
studied in \cite{Regular} (in which every element was encoded by its number
of lower successors).  In fact, if $F(\xi)$ denotes the latter encoding and
$\bar{\xi}$ is the permutation obtained from $\xi$ by replacing each element
$x_i$ by $|\xi|-x_i+1$ and then reversing it, then $F(\bar{\xi})$ is the
reverse of $E(\xi)$.  It therefore follows from Theorem 2 of \cite{Regular}
that

\begin{lemma}$E(X)$ is a regular set.
\end{lemma}

We can now confirm one of the claims we made when stating Theorem \ref{main}:
from the results of \cite{Regular} the ordinary generating function
of the sequence $(g_n)$, where $g_n$ is the number of permutations of length $n$,
is a rational function.  But 
to show that $\pi$ is eventually periodic and thus
complete the proof of the theorem we need
to study a deterministic finite automaton that accepts $E(X)$.  We denote
this automaton by $\M=(\Sigma,A,s_0,\tau,T)$ (the notation specifies,
respectively, the set of states, the alphabet, the initial state, the transition
function, and the set of final states).

\begin{lemma}\label{alphabeta}
$\{E(\sigma(i))\mid i=1,2,\ldots\}$ contains the
set 
\[\{\alpha\beta^j\mid j=0,1,2,\ldots\}\]
defined by the regular expression
$\alpha\beta^*$ for some non-empty words $\alpha,\beta\in A^\ast$.
\end{lemma}
\begin{proof}
We shall consider the sequence of states $t_i=\tau(s_0,E(\sigma(i))), i=1,2,\ldots$
and aim to show that it is periodic.
By definition and by Lemma \ref{uniqueness1} $\sigma(i+1)=\theta\phi$ where
\begin{enumerate}
\item[(a)]$\theta$ is a sequence order isomorphic to $\sigma(i)$
\item[(b)]$\phi=a_1a_2\cdots a_n$ satisfies $a_n<\max(\theta)$ and $\theta\phi\in X$
\item[(c)]$\phi$ is maximal with these properties.
\end{enumerate}
We shall express these conditions in terms of the automaton $\M$.  From (a),
(b) and the definition of $E$ we have
\[E(\sigma(i+1))=E(\sigma(i))b_1b_2\cdots b_n\]
where $w=b_1b_2\cdots b_n$ is a word in the alphabet $A$ and
$\tau(t_i,w)=t_{i+1}$.

By definition, $b_j$ is the number of terms of $\theta\phi$ up to and including
$a_j$ that exceed or equal $a_j$.  To capture the condition $a_n<\max(\theta)$
we need to define another sequence $c_0c_1c_2\cdots c_n$ where $c_0=0$
and $c_j$ is the
number of terms up to and including $a_j$ that exceed $\max(\theta)$; of
course, all such terms are among $a_1,a_2,\ldots,a_n$.

If $a_j<\max(\theta)$ then the terms enumerated by $b_j$ include 
$\max(\theta)$, $a_j$ and the $c_{j-1}$ terms above $\max(\theta)$; hence
$b_j>c_{j-1}+1$.  However, if $a_j>\max(\theta)$ then each of the $b_j$ terms that
exceed or equal $a_j$ is one of the terms that exceeds $\max(\theta)$; hence
$b_j\leq c_{j-1}+1$.  Furthermore, in the former case $c_j=c_{j-1}$ and in the latter
case $c_j=c_{j-1}+1$.  Thus $c_1c_2\cdots c_n$ is determined uniquely by
$b_1b_2\cdots b_n$, and $a_n<\max(\theta)$ if and only if $b_n>c_n$.

Putting all this together
$t_{i+1}$ is the unique state of $\M$ for which there exists a word
$w=b_1b_2\cdots b_n$ in the alphabet $A$ with the following properties:
\begin{enumerate}
\item $\tau(t_i,w)=t_{i+1}$;
\item if the sequence $(c_0,c_1,\ldots,c_{n})$ is defined by
$c_0=0\mbox{ and, for }j>0,$
\[c_{j}=\begin{cases}
c_{j-1}+1&\text{if $b_{j}\leq c_{j-1}+1$},\\
c_{j-1}&\text{if $b_{j}>c_{j-1}+1$}
\end{cases}
\]
then $b_n>c_{n}$;
\item $w$ has maximal length among all words satisfying these two conditions.
\end{enumerate}

But now note that the three conditions depend only on the $t_i$ and the automaton $\M$ and
not on $E(\sigma(i))$.  Therefore the sequence $t_1,t_2,\ldots$ is ultimately
periodic.  So, for some $P>0$ and $N$ we have $t_j=t_{j+P}$ for all $j\geq N$.

Let $\alpha=E(\sigma(N))$ and let $\beta$ be the unique word such that
$E(\sigma(N+P))=\alpha\beta$.  Then $E(\sigma(N+hP))=\alpha\beta^h$
and $\tau(s_0,\alpha\beta^h)=t_N$ for all $h\geq 0$.  This proves that
$\alpha\beta^*\subseteq\{E(\sigma(i))\mid i=1,2,\ldots\}$.
\end{proof}

We can now complete the proof of Theorem \ref{main}.  In the notation
of the previous lemma, let $|\alpha|=m$ and $|\beta|=n$.  Consider
the encoding of $\pi$ itself: $E(\pi)=e_1e_2\cdots$.  This is just
the limit of its prefixes $E(\sigma(1)),E(\sigma(2)),\ldots$.  It is also
the limit of $E(\sigma(N)), E(\sigma(N+P)), E(\sigma(N+2P)),\ldots$.  Hence,
by Lemma \ref{alphabeta}, $E(\pi)$ is ultimately periodic with 
$e_{j+P}=e_j$ for all $j\geq N$.

Consider an arbitrary $\pi(j)$ with $j\geq N$.  We have $\pi(j)=l_j+r_j+1$
where
\[l_j=|\{i\mid i\leq j\mbox{ and }\pi(i)<\pi(j)\}|\]
and
\[r_j=|\{i\mid i> j\mbox{ and }\pi(i)<\pi(j)\}|\]
Obviously, $l_j=j-e_j$.  The number $r_j$ can be obtained from $E(\pi)$ as 
follows.  Define two sequences $\nu^{(j)}=(n_0,n_1,\ldots)$ and 
$\theta^{(j)}=(h_0,h_1,\ldots)$ by $n_0=0$ and $h_0=e_j$ 
and
\[n_{i+1}=
\begin{cases}n_i+1&\text{if $e_{j+i+1}>h_i$}\\
n_i&\text{otherwise}
\end{cases}
\]
and
\[h_{i+1}=
\begin{cases}h_i&\text{if $e_{j+i+1}>h_i$}\\
h_i+1&\text{otherwise}
\end{cases}
\]
An easy inductive argument shows that $n_i$ is equal to the number of terms
from $\pi(j+1),\ldots,\pi(j+i)$ which are smaller than $\pi(j)$, while
$h_i-e_j$ is the number of terms from the same set which are greater
than $\pi(j)$.  In particular, $\nu^{(j)}$ eventually becomes constant with
value $r_j$.

Finally, note that $\nu^{(j)}$ depends only on $e_je_{j+1}e_{j+2}\cdots$ and
not on $e_1e_2\cdots e_{j-1}$.  Hence $\nu^{(j)}=\nu^{(j+P)}$ and $r_j=r_{j+P}$.
Therefore
\begin{eqnarray*}
\pi(j+P)&=&l_{j+P}+r_{j+P}+1\\
&=&j+P-e_{j+P}+r_{j+P}+1\\
&=&P+j-e_j+r_j+1\\
&=&P+l_j+r_j+1\\
&=&P+\pi(j)
\end{eqnarray*}
as required.

\section{Natural classes with infinite bases}\label{example}

Any natural class class that is not of the form stipulated by the conclusion 
of Theorem \ref{main} is, of course, not finitely based.  An example of such a class
is $Y=\Sub(\pi)$ where
\[
\pi = 3\ 2\ 5\ 1\ [7, 8]\ 4\ [10, 12]\ 6\ [14, 17]\ 9\ [19, 23]\ 13\ [25, 30]\ 18
\ [32, 38]\ 24\ldots\]
In this example $[a,b]$ stands for the segment $[a,a+1,\ldots,b]$.  It is clear by 
inspection that $\pi$ is not periodic and so, by the uniqueness conclusion of 
Theorem \ref{main}, $Y$ is not of periodic type.  We argue that it is not of the 
form $\Sub(\gamma)\oplus S$ where $S$ is sum-complete.  Suppose it were of
this form.  Consider the initial segments of $\pi$ ending with one of
$1,4,6,9,\ldots$ respectively.  These all define permutations 
$\xi_1,\xi_2,\xi_3,\xi_4\ldots$ of $Y$.
Every $\xi_i$ is indecomposable and has a unique embedding in $\pi$.  
From the indecomposability
those $\xi_i$ of length greater than $\gamma$ must be order isomorphic to permutations
of $S$; but, if $\xi_i\in S$, so also is $\xi_i\oplus\xi_i$ which
contradicts that it  is 
uniquely embeddable in $\pi$.

It would perhaps be tempting to suppose that when $\pi$ is periodic, the closed class $\Sub(\pi)$
is always finitely based.
This, however, is not the case, as our final example shows.

Let 
$X=\Sub(\pi)$, where 
\[
\pi= \underline{2\;3}\;\underline{5}\;1\;\underline{7\;8}\;4\;\underline{10}\;6\;\underline{12\;13}\;9\;\underline{15}\;11\cdots.
\]
Essentially, $\pi$ is an increasing oscillating sequence with every
other left maximal term replaced with an increasing pair (the
underlinings are intended to highlight this). 
Call these increasing pairs {\em twins}, and note that they are
the only pairs of terms of $\pi$ occurring in successive positions and having successive values.

We claim that each of the following permutations belongs to the basis of $X$:
\begin{eqnarray*}
\beta_1 &=&
\underline{2\;3}\;\underline{4\;5}\;1\\
\beta_2 &=&
\underline{2\;3}\;\underline{5}\;1\;\underline{7}\;4\;\underline{8\;9}\;6\\
&\vdots&\\
\beta_n &=&
\underline{2\;3}\;\underline{5}\;1\;\underline{7}\;4\;\underline{9}\;6\cdots \underline{4n-3}\; 4n-6\; \underline{4n-1}\; 4n-4\; \underline{4n\;4n+1}\;4n-2\\
&\vdots&
\end{eqnarray*}

\begin{figure}
\begin{center}
\begin{tabular}{ccc}
\begin{picture}(90,70)(0,0)
\thicklines
\put(0,0){\vector(1,0){70}}
\put(0,0){\vector(0,1){70}}
\put(3,6){\circle*{2}}
\put(6,9){\circle*{2}}
\put(12,15){\circle*{2}}
\put(15,3){\circle*{2}}
\put(18,21){\circle*{2}}
\put(21,24){\circle*{2}}
\put(24,12){\circle*{2}}
\put(27,30){\circle*{2}}
\put(30,18){\circle*{2}}
\put(33,36){\circle*{2}}
\put(36,39){\circle*{2}}
\put(39,27){\circle*{2}}
\put(42,45){\circle*{2}}
\put(45,33){\circle*{2}}
\put(48,51){\circle*{2}}
\put(51,54){\circle*{2}}
\put(54,42){\circle*{2}}
\put(57,60){\circle*{2}}
\put(65,40){\shortstack{\emph{etc}\\ \emph{ad} \\ \emph{inf}}}

\end{picture} & 

\begin{picture}(35,35)(0,0)
\thicklines
\put(0,0){\vector(1,0){35}}
\put(0,0){\vector(0,1){35}}
\put(3,6){\circle*{2}}
\put(6,9){\circle*{2}}
\put(12,15){\circle*{2}}
\put(15,3){\circle*{2}}
\put(18,21){\circle*{2}}
\put(21,12){\circle*{2}}
\put(24,27){\circle*{2}}
\put(27,30){\circle*{2}}
\put(30,18){\circle*{2}}

\end{picture} & 

\begin{picture}(45,45)(0,0)
\thicklines
\put(0,0){\vector(1,0){45}}
\put(0,0){\vector(0,1){45}}
\put(3,6){\circle*{2}}
\put(6,9){\circle*{2}}
\put(12,15){\circle*{2}}
\put(15,3){\circle*{2}}
\put(18,21){\circle*{2}}
\put(21,12){\circle*{2}}
\put(24,27){\circle*{2}}
\put(27,18){\circle*{2}}
\put(30,33){\circle*{2}}
\put(33,24){\circle*{2}}
\put(36,39){\circle*{2}}
\put(39,42){\circle*{2}}
\put(42,30){\circle*{2}}
\end{picture} 
\end{tabular}
\caption{On the left, an infinitely based periodic natural class.  On 
the right, two basis elements.}\label{inf.b.cycle.fig}
\label{fig4}
\end{center}
\end{figure}
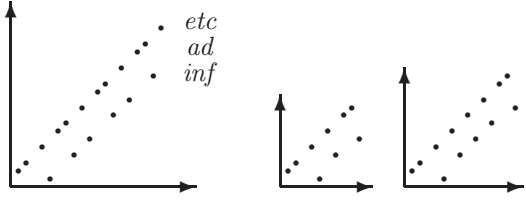

The permutation $\beta_n$ is obtained from an oscillating sequence with an even number of left maximal terms by replacing the first and last of these terms by increasing pairs (see Figure \ref{fig4}).
To show that $\beta_n\not\in X$ we can argue as follows.
Write $\pi$ as
$$
\pi=U_0U_1l_1U_2l_2U_3l_3\ldots,
$$
where $U_0,U_2,U_4,\ldots$ are the twins, $U_1,U_3,U_5,\ldots$ are
the remaining (unexpanded) left maxima, and $l_1,l_2,l_3,\ldots$ are
the remaining terms.
Suppose $\beta_n$ embeds into $\pi$.
The two twins of $\beta_n$ must correspond to two twins, say $U_{2p}$
and $U_{2q}$, of $\pi$.
Since $\beta_n$ is indecomposable, the remaining left maxima $5,7,9,\ldots$
of $\beta_n$ must map into $U_{2p+1},U_{2p+2},U_{2p+3},\ldots$ respectively.
The number of left maxima between the two twins of $\beta_n$ is even,
while the number of segments $U_{2p+1},\ldots,U_{2q-1}$ is odd,
a contradiction.

To complete the proof that $\beta_n$ is a basis permutation of $X$, we need to demonstrate that $\beta_n\setminus\{\beta_n(j)\}$, the permutation obtained 
by removing the $j$th term from $\beta_n$, belongs to $X$ for every
$j=1,\ldots,4n+1$.
If $j\not\in\{2,3,4n,4n+1\}$ the resulting permutation is decomposable, and can be embedded into $\pi$
by embedding each of its components and keeping them sufficiently apart.
If $j$ is one of $2,3,4n$ or $4n+1$ then one of the twins of $\beta_n$
becomes a singleton.
Suppose, for the sake of argument, that $j=4n+1$ (the other cases are
treated analogously).
Then we can embed $\beta_n\setminus\{\beta_n(4n+1)\}$
by mapping $2$ and $3$ onto $U_0$, all the other left maxima of $\beta_n$ into $U_1,U_2,\ldots$ respectively, and the remaining terms into 
$l_1,l_2,\ldots$ respectively.
Note that the parity problem which prevented us from embedding $\beta_n$ into
$\pi$ does not arise here, because the second twin of $\beta_n$ has become a singleton in $\beta_n\setminus\{\beta_n(4n+1)\}$,
and can therefore be mapped onto the singleton $U_{2n-1}$.

\end{document}